\documentstyle[12pt,epsfig]{article}
\newtheorem{thm}{Theora}[section]
\newtheorem{Theo}[thm]{Theorem}

\newtheorem{Lem}[thm]{Lemma}

\newtheorem{Def}[thm]{Definition}

\newenvironment{Proof}{\smallskip\noindent{\bf Proof.}\rm}
{\hfill $\Box$\medskip}

\newcommand{\eqref}[1]{(\ref{#1})}
\renewcommand\({\left(}
\renewcommand\){\right)}

\renewcommand\]{\right]}
\renewcommand\d{\partial}

\newcommand\R{{\bf R}}

\renewcommand\H{{\rm H}}

\renewcommand\L{{\rm L}}
\newcommand\la{\lambda}
\newcommand\ep{\varepsilon}

\newcommand\const{{\rm const}}

\makeatletter \@addtoreset{equation}{section}

\makeatother

\title{Systems-Disconjugacy of a Fourth-Order Differential Equation with a Middle Term}

\author{Jamel Ben Amara \\Facult\'e des Sciences de Bizerte, Tunisia\\
jamel.benamara@fsb.rnu.tn}

\date{}
\begin{document}
\maketitle
\begin{abstract}
Systems-conjugate points have been introduced and studied by John
Barrett \cite{barr} in relation with the self-adjoint fourth order
differential equation $$(r(x)y'')''-(q(x)y')'= p(x)y,$$ where
$r(x)>0$, $p(x)>0$ and $q\equiv0$. In this paper we extend some of
his results to more general cases, when $q(x)$ is free of any sign
restrictions.
\end{abstract}



\section{Introduction}
This paper shall be concerned with the fourth-order differential
equation
\begin{equation}\label{eq:mainq}
(r(x)y'')''-(q(x)y')'= p(x)y,
\end{equation}
where $r(x)>0$, $p(x)>0$ and $q(x)$ are continuous functions on
$[a,\infty )$, $a\geq0$.

\begin{Def}\label{eq:def}
 The systems-conjugate
point of $a$, which is denoted by $\hat\eta_1(a)$, is defined as
the smallest number $b\in(a,\infty )$ for which the two point
boundary conditions
\begin{equation}\label{eq:bc1}
y(a) = y_1(a)=y(b) = y_1(b)=0
\end{equation}
($y_1(x)=r(x)y''$) are satisfied by a nontrivial solution of
equation \eqref{eq:mainq}. \\
Similarly, the systems-focal point of $a$, which is denoted by
$\hat\mu_1(a)$, is defined as the smallest number $b\in(a,\infty
)$ for which the two point boundary conditions
\begin{equation}\label{eq:bc2}
y(a) = y_1(a)=y'(b) = Ty(b)=0
\end{equation}
 ($Ty(x)= (p(x)y'')' - q(x)y'$) are satisfied by
 a nontrivial solution of equation \eqref{eq:mainq}.
\end{Def}
The notation $y_1(x)$ and $Ty(x)$ will be used throughout the paper.\\
\indent The systems-conjugate point and systems-focal point were
first defined and studied by Barrett \cite{barr,barrAd} with respect
to equation \eqref{eq:mainq}, for $r(x)>0$, $p(x)>0$ and $q\equiv0$.
In his work, he showed that $\hat\eta_1(a)$ exists, if and only if
$\hat\mu_1(a)$ exists, and $a<\hat\mu_1(a)<\hat\eta_1(a)$, without
further conditions on $r(x)$ and $p(x)$. Later on, using a Morse
system-formulation \cite{Morse}, Atkinson \cite[Chap. 10.6]{Atk}
extended a part of Barrett's result to the case $q(x)\geq0$ (i.e.,
if $\hat\eta_1$ exists then $\hat\mu_1(a)$ exists and
$a<\hat\mu_1(a)<\hat\eta_1(a)$). Cheng \cite{cheng} also studied the
existence and the relation between $\hat\mu_1(a)$ and
$\hat\eta_1(a)$ for a system of two second-order differential
equations; in particular, he gave a physical interpretation of the
numbers $\hat\eta_1$ and $\hat\mu_1$. At the end of this work, he
applied his results to equation \eqref{eq:mainq} for $q(x)\leq0$ and
 the additional condition $p -q''/2 +{q^2}/4r
>0$. Note that the systems-focal point studied in \cite{cheng} do not
coincide with that defined above for \eqref{eq:mainq} only for
$q\equiv \const$.\\
 \indent The main goal of the present paper is to establish Barrett's
result related to equation \eqref{eq:mainq} with some relaxation of
the sign of $q(x)$. Furthermore, in Sections 3 and 4 we establish a
comparison theorem for $\hat\mu_1(a)$, and we show, without further
restrictions on $r$, $p$ and $q$, that if $\hat\mu_1(a)$ exists then
it is realized by a positive increasing solution. These results are
analogous to those obtained by Barrett \cite{barrd} for the focal
point $\mu_1(a)$ related to equation \eqref{eq:mainq} and the
boundary conditions $y(a) = y'(a)=y_1(b) = Ty(b)=0$. However, here
we use a different approach, which is essentially based on the
Leighton-Nehari transformation \cite{LN} and the properties of the
Rayleigh quotients. Finally, in Section~5 we establish two criteria
for the existence of $\hat\eta_1(a)$. Similar results were given in
\cite{barr} and \cite{cheng} for $q(x)\equiv0$ and $q(x)\leq0$,
respectively.

\section {Relation between $\hat\eta_1$ and $\hat\mu_1$}
The main result of this section is the following

\begin{Theo}\label{le:main}
$1)$ If the first systems-conjugate point $\hat\eta_1(a)$ exits
and
\begin{equation}\label{eq:quadform}
I(w,a,b)=\int_a^b[ r(w')^2 + qw^2] >0
\end{equation}
for each $b>a$ and each nontrivial admissible function $w\in
W_2^1[a,b]$ (where $W_2^1[a,b]$ is the Sobolev function space
having a generalized first derivative in $\L_{2}[a, b]$), then the
first systems-focal point $\hat\mu_1(a)$ exists and
\begin{equation}\label{eq:ineg}
a< \hat\mu_1(a)<\hat\eta_1(a).
\end{equation}
$2)$ If the number $\hat\mu_1(a)$ exists and $\int^\infty
q(t)=-\infty$, then $\hat\eta_1(a)$ exits. If in addition the
condition \eqref{eq:quadform} is satisfied, then \eqref{eq:ineg}
holds.
\end{Theo}

 Before proving this theorem we need some preliminaries.
 It is known that any solution of equation \eqref{eq:mainq} which
satisfies the initial condition $y(a) = y_1(a)=0$ may be expressed
as a linear combination of $u(x)$ and $v(x)$ which are the
fundamental solutions of \eqref{eq:mainq} whose initial conditions
are
\begin{equation}\label{eq:initu}
u(a) = u_1(a) = Tu(a)=0,\quad u'(a)=1,
\end{equation}
\begin{equation}\label{eq:initv}
v(a) = v'(a)=v_1(a)=0,\quad  Tv(a)=1.
\end{equation}
We introduce the following subwronskians:
\begin{equation}\label{eq:wrons}
r\hat\sigma'= uv_1 -vu_1, \quad \hat\tau'= u'Tv - v'Tu,
\end{equation}
and
\begin{equation}\label{eq:wrons1}
\hat\sigma= uv' -vu', \quad \hat\tau= uTv - vTu,\quad \hat\rho=
u_1Tv - v_1Tu.
\end{equation}
It is easy to see that $\hat\eta_1$ and $\hat\mu_1$ are the first
zeros on $(a,\infty)$ of the subwronskians $\hat\sigma'$ and
$\hat\tau'$, respectively. The following identities involving the
above subwronskians are useful and easily verified. Similar ones
have been stated in \cite{barrd} for \eqref{eq:mainq} with Dirichlet
boundary conditions at $x=a$ ($y(a) = y'(a)=0$).

\begin{equation}\label{eq:ident1}
r\hat\sigma'\hat\tau'=\hat\tau^2 + \hat\rho\hat\sigma
\end{equation}

\begin{equation}\label{eq:ident2}
\hat\tau''=\frac{\hat\rho}{r}- p\hat\sigma,\quad  (r\hat\sigma')'
= 2\hat\tau + q\hat\sigma.\end{equation}

Note also, the initial conditions
\begin{equation}\label{eq:initc1}
\hat\tau(a)=0,\quad \hat\tau'(a)=1,
\end{equation}
\begin{equation}\label{eq:initc2}
\hat\sigma(a)=\hat\sigma'(a)=(r\hat\sigma')'(a)=0,
(r\hat\sigma')''(a)=2,
\end{equation}
\begin{equation}\label{eq:initc3}
\hat\rho(a)=0, \quad \hat\rho'(a)=q(a),
\end{equation}
insure that $\hat\sigma$, $\hat\sigma'$, $\hat\tau$ and $\hat\tau'$
are all positive in a right-hand neighborhood of
$x=a$. \\
Throughout our discussion we will use the following
transformation given by Leighton-Nehari \cite{LN} for removing the
middle term $(qy')'$ from equation \eqref{eq:mainq}. However, this
transformation can not be used in a straightforward way, since as
will be seen below, it changes the form of the initial conditions
\eqref{eq:initu} and the subwronskians $\hat\sigma'$, $\hat\rho$.\\
Let us denote by $h$ a positive solution on the interval $[a, b]$ of
the second-order equation
\begin{equation}\label{eq:disf}
(py')'- qy=0.
\end{equation}
Hence, the following substitution \cite[Theorem 12.1]{LN}
\begin{equation}\label{eq:chv}
t(x):={\int_0^x h(s)ds}
\end{equation}
 transform equation \eqref{eq:mainq} into
\begin{equation}\label{eq:transeq}
{\(r h^{3}(t)\ddot y\)}^{..} =  h^{-1}p(t) y,\\
\end{equation}
where $p(x),h(x),r(x),y(x)$ are taken as functions of $t$ and
$^\cdot:=\frac{d}{dt}$. Therefore, if $y$ is a nontrivial solution
of \eqref{eq:mainq}, then $\tilde y(t)\equiv y(x(t))$ is a
nontrivial solution of \eqref{eq:transeq}. Thus, we have the
relations:
\begin{equation}\label{eq:re1}
\dot {\tilde y} =y'h^{-1},\quad h^{3}{\ddot {\tilde y}}=hy''- y'h',
\end{equation}
\begin{equation}\label{eq:re2}
{(\tilde r\tilde h^{3}\ddot {\tilde y})}^{.}=(ry'')'- qy'.
\end{equation}
In what follows, for each of the quantities involving
\eqref{eq:transeq}, the same notations as for \eqref{eq:mainq} will
be used with the addition of the superscript ``$\sim$". Let
$\tilde\sigma$, $\dot{\tilde\sigma}$, ${\tilde\tau}$,
$\dot{\tilde\tau}$ and $\tilde\rho$ denote the subwronskians
associated with equation \eqref{eq:transeq} and the fundamental
solutions $\tilde u$, $\tilde v$ satisfying the initial conditions
\begin{equation}\label{eq:tilde u}
\tilde {u}(0) = \tilde h(0) \ddot {\tilde u}(0)+\dot {\tilde h}(0)
= {(\tilde r\tilde h^{3}\ddot {\tilde u})}^{.}(0)=0,\; \dot{\tilde u}(0)=1,
\end{equation}
\begin{equation}\label{eq:tilde v}
\tilde {v}(0) = \dot {\tilde v}(0)= \ddot {\tilde v}(0) = 0,\;
{(\tilde r\tilde h^{3}\ddot {\tilde v})}^{.}(0)=1.
\end{equation}
The relations between these subwronskians and those of equation
\eqref{eq:mainq} are expressed as follows:
\begin{equation}\label{eq:relat1}
\hat\sigma(x)= \tilde h(t)\tilde \sigma(t),\quad \hat\tau(x)= \tilde
{\tau}(t), \quad \hat\tau'(x)=\tilde h(t)\dot {\tilde{\tau}}(t),
\end{equation}
\begin{equation}\label{eq:relat2}
\hat\sigma'(x)={\tilde h^2(t)\dot{\tilde\sigma}(t)} + \tilde
 h\dot{\tilde h}\tilde\sigma(t),
\end{equation}

\begin{Lem}\label{le:rho}
1) If $\hat \mu_1(a)$ exists, then $\hat\rho(\hat\mu_1(a))<0$.\\
2) Let $\hat\xi_1(a), \hat\xi_2(a)\cdots$ denote the zeros of the
subwronskian $\hat \rho$ defined by \eqref{eq:wrons}. If
$\hat\mu_1(a)$ exists and $\hat\rho$ has a first zero $\hat
\xi_i(a)$ ($i\in\{2, 3, \cdots\}$) beyond $\hat\mu_1(a)$, then
$\hat\sigma'$ has a zero $\hat \eta_1(a)$ in $(a, \hat \xi_i(a)]$.
\end{Lem}
\begin{Proof}
$1)$ If $\hat \mu_1(a)$ exists, then $\hat\sigma>0$ on $(a, \hat
\mu_1(a)]$. In fact, suppose that $\hat\sigma$ has a zero $s_0\in(a,
\hat \mu_1(a))$ which is the closest to $a$. From the initial
conditions \eqref{eq:initc1}-\eqref{eq:initc2}, we have
$\hat\tau'>0$ and $\hat\sigma>0$ in a right-hand neighborhood of
$x=a$, and hence, $\hat\tau'(s_0)>0$ and $\hat\sigma'(s_0)<0$. On
the other hand, by \eqref{eq:ident1},
$\hat\tau'\hat\sigma'(s_0)\geq0$, which is a contradiction. If $s_0=
\hat \mu_1(a)$, then again by \eqref{eq:ident1}, $\hat\tau(s_0)=0$.
Thus, from Rolle's theorem and the initial conditions
\eqref{eq:initc1}, there exists
a zero of $\hat\tau'$ less than $\hat \mu_1(a)$, which is a contradiction.\\
Since $\hat\sigma(\hat \mu_1(a))>0$ then by \eqref{eq:ident1}, we
have $\hat\rho\hat\sigma(\hat \mu_1(a))\leq0$. If $\hat\rho(\hat
\mu_1(a))=0$, then $\hat\tau(\hat \mu_1(a))=0$, and
as before, this is not possible.\\
$2)$ Suppose that $\hat\rho$ has a first zero $\hat\xi_i$ beyond
 $\hat \mu_1(a)$ (i.e., the first in $(\hat\mu_1(a), \infty)$).
By \eqref{eq:ident1}, we have
$\hat\tau'\hat\sigma'(\hat\xi_i(a))\geq0$.
 If $\hat\tau'(\hat\xi_i(a))<0$, then $\hat\sigma'(\hat\xi_i(a))\geq0$, and hence, from the
initial conditions \eqref{eq:initc2}, $\hat\sigma'$ has a zero $\hat
\eta_1(a)$ in $(a, \hat \xi_i(a)]$. If $\hat\tau'(\hat
\xi_i(a))\geq0$, then $\hat \mu_2(a)$ exists and
$a<\hat\mu_2\leq\hat \xi_i(a)$. According to Lemma \ref{le:sigm},
$\hat\sigma$ has a zero in the interval $(\hat \mu_1(a), \hat
\mu_2(a)]$. Thus, by Rolle's theorem, $\hat\eta_1(a)$ exists and
$a<\hat\eta_1(a)\leq\hat \xi_i(a)$. The lemma is proved.
\end{Proof}
\begin{Lem}\label{le:sigm}
If $\hat \mu_1(a)$ and $\hat\mu_2(a)$ (the second zero of $\tau'$)
both exist, then $\hat\sigma$ has a zero in the interval $(\hat
\mu_1(a), \hat \mu_2(a)]$.
\end{Lem}
\begin{Proof}
By Lemma \ref{le:rho} and its proof, we have $\hat\sigma(\hat
\mu_1(a))>0$ and $\hat\rho(\hat \mu_1(a))<0$. Thus, $\hat\tau''(\hat
\mu_1(a))=(\frac{\hat\rho}{r}- p\hat\sigma)(\hat \mu_1(a))<0$, which
implies the simplicity of $\hat \mu_1(a)$, and hence, $\hat
\mu_1(a)<\hat\mu_2(a)$. Suppose $\hat\sigma>0$ on $(\hat \mu_1(a),
\hat \mu_2(a)]$. Since $\hat\tau(\hat \mu_1(a))>0$, then by using
the identity \eqref{eq:ident1}, we obtain
$$
\(\frac{\hat\tau'}{\hat\sigma}\)' =-p -
\frac{1}{r}\(\frac{\hat\tau}{\hat\sigma}\)^2<0.
$$
Integration of this expression yields
$$
\int_{\hat \mu_1}^{\hat \mu_2} p +
\frac{1}{r}\(\frac{\hat\tau}{\hat\sigma}\)^2 dx=0,
$$
which is a contradiction, and so $\hat\sigma$ vanishes in $(\hat
\mu_1(a), \hat \mu_2(a)]$.
\end{Proof}

\begin{Proof}{\bf {of Theorem \ref{le:main}}}\\

$1)$ Let $h$ be the solution of equation \eqref{eq:disf} which
satisfies the initial conditions
\begin{equation}\label{eq:initial}
y'(a)=0,\quad y(a)=1.
\end{equation}
If condition \eqref{eq:quadform} holds, then all the eigenvalues of
the problem determined by equation \eqref{eq:disf} and the boundary
conditions $y'(a)=y'(b)=0$ (for each $b>a$) are positive, and hence,
$h(x)>0$ on $[a, \infty)$. Furthermore, since
$I(1,a,a+\ep)=\int_a^{a+\ep} q >0$ for sufficiently small $\ep>0$,
$q(x)\geq0$ in a right-neighborhood of $x=a$. Thus, $h'(x)>0$ on
$[a, \infty)$. Therefore, the change of variables $t(x):=\int_0^x
h(s)ds$ is valid to transform equation \eqref{eq:mainq} into
\eqref{eq:transeq}. Let $\tilde\eta_1(0)$ and $\tilde\mu_1(0)$
denote, respectively, the first systems-conjugate point and the
first systems-focal point associated with equation
\eqref{eq:transeq}; i.e., the first zeros of the subwronskians
$\dot{\tilde\sigma}$ and $\dot{\tilde\tau}$, respectively. As noted
before, these subwronskians are obtained from the original ones via
the above change of variables, and the relations between them are
expressed by \eqref{eq:relat1}-\eqref{eq:relat2}. Note that also the
initial conditions
\begin{equation}\label{eq:initc3'}
\tilde\tau(0)=0,\quad \dot{\tilde\tau}(0)=1,
\end{equation}
\begin{equation}\label{eq:initc4}
{\tilde\sigma}(0)=\dot{\tilde\sigma}(0)=(\tilde{r}\tilde{h}^3
\dot{\tilde\sigma})^{.}(0)=0, (\tilde{r}\tilde{h}^3
\dot{\tilde\sigma})^{..}(0)=2,
\end{equation}
imply that ${\tilde\sigma}$, $\dot{\tilde\sigma}$ and
$\dot{\tilde\tau}$ are positive in a right-hand neighborhood of
$t=0$.\\
Suppose $\hat\eta_1(a)$ exists. By \eqref{eq:initial}, together with
the relation \eqref{eq:relat2}, we have
$\dot{\tilde\sigma}(\int_a^{\hat\eta_1(a)}h)<~0$. Hence,
$\tilde\eta_1(0)$ exists for \eqref{eq:transeq}. According to
\cite[Theorem 1.1]{barr}, which is applied to equation
\eqref{eq:transeq}, it follows that $\tilde\mu_1(0)$ exists, and
\begin{equation}\label{eq:ineq4}
0< \tilde\mu_1(0)< \tilde\eta_1(0).
\end{equation}
Therefore, from the last relation of \eqref{eq:relat1},
$\hat\mu_1(a)$ also exists, and \eqref{eq:ineg} holds.\\
$2)$ Assume that $\int^{\infty}q=-\infty$, and suppose that
$\hat\mu_1(a)$ exists, but $\hat\sigma'>0$ on $(a, \infty)$. In view
of Lemmas \ref{le:rho} (second statement) and \ref{le:sigm}, if
$\hat\xi_i(a)$ (the first zero of $\hat\rho$ beyond $\hat\mu_1(a)$)
or $\hat\mu_2(a)$ exists, then $\hat\eta_1(a)$ exists. On the other
hand, by the first statement of Lemma \ref{le:rho},
$\hat\rho(\hat\mu_1(a))<0$. Therefore, if $\hat\xi_i(a)$ and
$\hat\mu_2(a)$ do not exist then we have
$k(x)=-\frac{\hat\rho}{\hat\tau'}<0$ on $(\hat\mu_1(a), \infty)$,
and
$$
k'(x) =p(\frac{\hat\tau}{\hat\tau'})^2 - q + \frac{1}{r}k^2\geq0
\qquad\mbox{on}\qquad(a, \infty).
$$
Integrating this expression, and taking into account the assumption
that $\int^{\infty}q=-\infty$, it follows that
$k(x)\rightarrow+\infty$ as $x\rightarrow+\infty$. This is a
contradiction, and so $\hat\mu_1(a)$ exists. If in addition,
\eqref{eq:quadform} holds, then from the first statement of the
theorem, we have $a< \hat\mu_1(a)<\hat\eta_1(a)$. The theorem is
proved.
\end{Proof}

\section { Wirtinger
inequality and comparison theorem for $\hat\eta_1(a)$ }

The following theorem establish the relation between the existence
of $\hat\eta_1(a)$ and the sign of the quadratic form associated
with \eqref{eq:mainq}. This relation is known as a Wirtinger-type
inequality \cite{cole}. Note that the method of Cole used in
\cite{cole} and also in \cite{barrd} (for a Wirtinger inequality
related to the focal point $\mu_1(a)$) cannot be applied here.

\begin{Theo}\label{le:quadrat}
If $\hat\eta_1(a)$ does not exist for \eqref{eq:mainq}, then for
each $b\in\(a, +\infty\)$ and each nontrivial admissible function
$w(x)$ on $[a, b]$ (i.e., $w(x)\in C^{1}[a, b]$, $w'$ is absolutely
continuous and $w''\in L_{2}[a, b]$) for which $w(a)=w'(b)=0$, we
have
$$
I[w, a, b]={\int_a^b r(w'')^2 + q(w')^2 - pw^2 dx}>0.
$$
\end{Theo}

For the proof of this theorem, we need some preliminarily results.
We introduce the following equation similar to \eqref{eq:mainq},
but depends on a parameter $\la\in\R$.
\begin{equation}\label{eq:mainqq}
(r(x)y'')''-(q(x)y')'=\la p(x)y.
\end{equation}
Let $\bar\eta_1(a)$ denote the first conjugate point of $a$ with
respect to equation \eqref{eq:disf}; i.e., the smallest number
$b\in(a,\infty )$ for which the boundary conditions $y'(a)=y(b)=0$
are satisfied by a nontrivial solution.

\begin{Lem}\label{le:first}
Let $\la_1(b)$ be the first eigenvalue  of Problem
\eqref{eq:mainqq}-\eqref{eq:bc2}, and assume that $\bar\eta_1(a)$
exists. If $\la_1(b)>0$, then $b< \bar\eta_1(a)$.
\end{Lem}

The proof of this lemma is based on the following result on the
monotonicity of the eigenvalues of Sturm-Liouville problem. To the
best of my knowledge, this property is known only for $q\geq0$
(e.g., \cite{courant}).

\begin{Lem}\label{le:first-second}
The eigenvalues $\rho_k(b)$ of the second-order boundary problem
\begin{equation}\label{eq:second-eigen}
-(r(x)y')'+ q(x)y=\rho y,
\end{equation}
\begin{equation}\label{eq:initial2}
y'(a)=0,\quad y(b)=0
\end{equation}
 decrease as $b$ increases.
\end{Lem}

\begin{Proof}
Let $$ F(x,\rho) = \frac{y(x,\rho)}{ry'(x,\rho)},
$$
where $y(x, \rho)$ is a nontrivial solution of Problem
\eqref{eq:second-eigen}-\eqref{eq:initial2}. Obviously, for fixed
$\rho$, the zeros and poles of $F(x,\rho)$ do not coincide unless
$y(x, \rho)\equiv0$. If $y(b,\rho_k(b))=0$, then $F(b,\rho_k(b))=0$
and
\begin{equation}\label{eq:F}
\frac{\d F(x,\rho_k(b))}{\d x}_{|x=b} = 1/r(b)>0.
\end{equation}
On the other hand, for fixed $x=b$, $F(x,\rho)$ is a finite-order
meromorphic function of $\rho$, and
\begin{eqnarray}\label{eq:ddF}
\frac{\d{F(b,\rho)}}{\d{\rho}}_{|\rho=\rho_k(b)}=
y^{-2}(b,\rho_k(b))\int_a^b p(x) y^2(x,\rho_k(b))dx >0
\end{eqnarray}
(e.g., see, \cite[Chap.6]{Atk}). From the implicit-function
theorem, together with \eqref{eq:F}-\eqref{eq:ddF}, we obtain
$$
\rho'_k(b)=-\frac{\frac{\d F(x,\rho_k(b))}{\d
x}_{|x=b}}{\frac{\d{F(b,\rho)}}{\d{\rho}}_{|\rho=\rho_k(b)}}<0,
$$
and this completes the proof of the lemma.
\end{Proof}

\begin{Proof}{\bf {of Lemma \ref{le:first}}}\\
Suppose $\la_1(b)>0$, but $b\geq \bar\eta_1(a)$. In this case, the
mini-max principle yields:
$$
\la_1(b)= \min_{w\in\H} \frac{I(w)}{\int_a^b p(w)^2 dx}>0,
$$
where $I(w)={\int_a^b \[r(w'')^2 + q(w')^2\]dx}$, and $\H$ is a set
of nontrivial admissible functions $w$ (i.e., $w(x)\in C^{1}[a, b]$,
$w'$ is absolutely continuous and $w''\in L_{2}[a, b]$) for which
$w(a)=w'(b)=0$. On the other hand, by Lemma \ref{le:first-second},
$\rho_1(b)\leq0$, and hence, the corresponding eigenfunction $v(x)$
satisfies the inequality
$$
{\int_a^b \[r(v')^2 + q(v)^2\] dx}\leq0.
$$
Let $y(x)=\int_a^x vdx$. Then $y(a)=0$, $y'(b)=0$ and $ {\int_a^b
\[r(y'')^2 + q(y')^2\] dx}\leq0$, which is a contradiction. The lemma is
proved.
\end{Proof}

The conclusion in the second part of the following lemma is similar
to that of Greenberg \cite{green} stated for the first eigenvalue of
the problem determined by equation \eqref{eq:mainqq} and the
Dirichlet boundary conditions $y(a) = y'(a)=y(b) = y'(b)=0$.

\begin{Lem}\label{le:infinity}
The first eigenvalue $\la_1(b)$ of Problem
\eqref{eq:mainqq}-\eqref{eq:bc2} is simple. Furthermore, if
$b\rightarrow +a$ then $\la_1(b)\rightarrow +\infty$.
\end{Lem}

\begin{Proof}
By Lemma \ref{le:first}, if $\la_1(b)>0$ then $b<\bar\eta_1(a)$.
Therefore, the solution $h$ of the initial-value problem
\eqref{eq:disf}-\eqref{eq:initial} is positive on the interval
$[a,b]$, and hence, it is possible to transform equation
\eqref{eq:mainqq} (with $\la=\la_1(b)$) into

\begin{equation}\label{eq:transeqq}
{\(r h^{3}(t)\ddot y\)}^{..} = \la h^{-1}p(t) y
\end{equation}
(with $\la=\la_1(b)$), and the boundary conditions \eqref{eq:bc2}
into
\begin{equation}\label{eq:bc4}
\tilde {y}(0) = \ddot {\tilde v}(0) =\dot {\tilde v}(\tilde b)=
{(\tilde r\tilde h^{3}\ddot {\tilde v})}^{.}(\tilde b)=0,
\end{equation}
where $\tilde b=\int_a^b h dx$. Obviously, if $\la=\la_1(b)$ is a
multiple eigenvalue of Problem \eqref{eq:mainqq}-\eqref{eq:bc2},
then it is so for Problem \eqref{eq:transeqq}-\eqref{eq:bc4}. But,
this is not possible since all the eigenvalues of this problem are
simple (e.g., see \cite{BK}).\\
Let $ b_0>a$. For each $b\leq b_0$, consider the quadratic form
$$I(y)={\int_a^b \[r(y'')^2 + q(y')^2\] dx},$$ defined on the set of
all nontrivial admissible functions $y$ (i.e., $y(x)\in C^{1}[a,
b]$, $y'$ is absolutely continuous and $y''\in L_{2}[a, b]$) for
which $y(a)=y'(b)=0$. For such $y$, we have the following
expressions, which follows from the Cauchy-Schwarz inequality:
$$
\int_a^b (y)^2 dx\leq (b-a)\int_a^b (y')^2 dx,
$$
and
$$
\int_a^b (y')^2 dx\leq (b-a)\int_a^b (y'')^2 dx.
$$
 Therefore,
$$
I(y)\geq  \frac{r^* \int_a^b (y)^2 dx}{(b-a)^2} + \frac{q^*
\int_a^b (y)^2 dx}{(b-a)},
$$
where, $f^*=\min_{x\in[a,b_0]}f(x)$. Thus,
$$
 \frac{I(y)}{\int_a^b p(y)^2 dx}\geq \frac{1}{p_*}\(\frac{r^*}{(b-a)^2} +
 \frac{q^*}{(b-a)}\),
$$
where, $p_*=\max_{x\in[a,b_0]}p(x)$. The mini-max principle
implies
$$
\la_1(b)\geq \frac{1}{p_*}\(\frac{r^*}{(b-a)^2} +
 \frac{q^*}{(b-a)}\),
$$
and hence, $\lim_{b\rightarrow a}\la_1(b)=+\infty$.
\end{Proof}

\begin{Proof}{\bf {of Theorem \ref{le:quadrat}}}
In view of Lemma \ref{le:infinity}, $\la_1(b)\rightarrow +\infty$ if
$b\rightarrow +a$ (recall that $\la=\la_1(b)$ denotes the smallest
eigenvalue of Problem \eqref{eq:mainqq}-\eqref{eq:bc2}). Thus, there
exists $b>a$ such that $\la_1(b)>1$. Let $\hat\tau'(\la, x)$ denotes
the subwronskian defined by \eqref{eq:wrons} related to equation
\eqref{eq:mainqq}. It is easily remarked that, for fixed $x=b$, the
zeros of the function $\hat\tau'(\la, x)$
 and the eigenvalues of Problem \eqref{eq:mainqq}-\eqref{eq:bc2},
together with their multiplicities, coincide. In particular, the
simplicity of $\la_1(b)$ (see Lemma \ref{le:infinity}) yields
$$
  \hat\tau'(\la_1(b), b)=0,\quad \frac{\d \hat\tau'}{\d {\la}}
 (\la, b)_{|\la=\la_1(b)}\ne0.
$$
It then follows from the implicit-function theorem that $\la_1(b)$
is a continuous function of $b\in(a,\infty)$. Therefore, as $b$
varies along the interval $(a,\infty)$, $\la_1(b)$ can not pass
through the value $\la=1$, since otherwise, we have for some $b>a$,
$\hat\eta_1(a)=b$ exists for \eqref{eq:mainq}, and this is in
contradiction to the hypothesis of the theorem. Hence, $\la_1(b)>1$
for all $b\in\(a, +\infty\)$, and so, for every nontrivial
admissible function $w$ for which $w(a)=w'(b)=0$, we obtain
$$
\int_a^b r(w'')^2 + q(w')^2 dx > {\int_a^b pw^2 dx}.
$$
The theorem is proved.
\end{Proof}

We now establish a comparison theorem for
 $\hat\mu_1(a)$.

\begin{Theo}\label{le:comp}
Let $r_0(x)>0$, $p_0(x))>0$ and $q_0(x)$ be continuous functions on
$[a, \infty)$, such that
\begin{equation}\label{eq:coeff}
r\leq r_0,\qquad p_0 \leq p, \qquad q_0 \geq q,
\end{equation}
and there exists the first systems-focal point, say
$\hat\mu^0_1(a)$, for the equation
\begin{equation}\label{eq:nequat}
(r_0(x)y'')''-(q_0(x)y')'= p_0(x)y.
\end{equation}
Then $\hat\mu_1(a)$ exists for the original equation
\eqref{eq:mainq} and
$$
a< \hat\eta_1(a)\leq\hat\eta^0_1(a).
$$
\end{Theo}
\begin{Proof}

Suppose that $\hat\eta^0_1(a)$ exists but $\hat\sigma'>0$ on $(a,
\hat\eta^0_1(a)]$. Let $y_0$ be the corresponding eigenfunction;
then Theorem \ref{le:quadrat} yields

$$
I[y_{0}, a, \hat\eta^0_1(a)]={\int_a^{\hat\eta^0_1(a)} r(y''_{0} )^2
+ q(y'_{0})^2 - p(y_{0} )^2 dx}>0
$$

and

$$
I^{0}[y_{0}, a, \hat\eta^0_1(a)]={\int_a^{\hat\eta^0_1(a)}
r_0(y''_{0} )^2 + q_0(y'_{0})^2 - p_0(y_{0} )^2 dx}=0.
$$

Subtracting these two expressions and taking into account
\eqref{eq:coeff}, we obtain
$$
0\leq{\int_a^{\hat\eta^0_1(a)} (r_0 -r)(y''_{0} )^2 + (q_0
-q)(y'_{0})^2 +(p- p_0)(y_{0} )^2 dx}<0.
$$
This contradiction shows that there exists
$\hat\eta_1(a)\leq\hat\eta^0_1(a)$.
\end{Proof}

\section {Oscillation of the eigenfunction associated to $\hat\mu_1(a)$}

\begin{Theo}\label{le:oscill}
If $\hat\mu_1(a)$ exists, then it is realized by an unique
eigenfunction $y_{\hat\mu_1}$ up to a multiplicative constant. It
has the properties
 $$ y_{\hat\mu_1}>0,\quad  y'_{\hat\mu_1}>0,\quad
T(y_{\hat\mu_1})<0\quad on\quad (a, \hat\mu_1).$$ Also, if
$q\leq0$ on $[a, \hat\mu_1)$, then $y''_{\hat\mu_1}<0$ on $(a,
\hat\mu_1)$.
\end{Theo}

The following lemma establishes the relation between $\bar\eta_1$
(defined in Section 3) and $\hat\mu_1$.

\begin{Lem}\label{le:barrett1}
If $\bar\eta_1(a)$ exists, then $\hat\eta_1(a)$ exists, and
\begin{equation}\label{eq:int1}
a < \hat \mu_1(a) \leq \bar \eta_1(a),
\end{equation}
with equality if, and only if, $p(x)\equiv0$ on $[a, \bar
\eta_1(a)]$.
\end{Lem}
\begin{Proof}
It is easily seen that if $p(x)\equiv0$, then $\hat \mu_1(a)=\bar
\eta_1(a)$. Therefore, the conclusion of the lemma follows from
Theorem \ref{le:comp}.
\end{Proof}

For the proof of Theorem \ref{le:oscill} we need the following two
lemmas.

\begin{Lem}\label{le:LN} {\rm ([LN, Lemma 2.1])}
Let $y$ be a nontrivial solution of the differential equation
\eqref{eq:mainq} for $q\equiv0$. If $y, y', y''$ and $Ty$ are
nonnegative at $x=a$ (but not all zero), then they are positive for
all $x>a$. If $y, -y', y''$ and $-Ty$ are nonnegative at $x=a$ (but
not all zero), then they are positive for all $x<a$.
\end{Lem}

\begin{Lem}\label{le:fund}
Let $u$ and $v$ be two fundamental solutions of \eqref{eq:mainq}
defined by \eqref{eq:initu} and  \eqref{eq:initv}, respectively.
Then:
\begin{equation}\label{eq:oscil1}
u>0,\quad u'>0,\quad Tu>0 \quad \mbox{on} \quad(a, \hat\mu_1]. \\
\end{equation}
\begin{equation}\label{eq:oscil2}
v>0,\quad v'>0,\quad Tv>0 \quad \mbox{on} \quad(a, \hat\mu_1].
\end{equation}
If, in addition; \eqref{eq:quadform} holds, then $u''>0$ and $v''>0$
on $(a, \hat\mu_1]$.
\end{Lem}
\begin{Proof}
In view of Lemma \ref{le:barrett1}, we have
$a<\hat\eta_1(a)<\bar\eta_1(a)$. In this case, from the definition
of $\bar\eta_1(a)$, the solution $h$ of \eqref{eq:disf} satisfying
the initial conditions $h'(a)=0$, $h(a)=1$, is positive on $[a,
\hat\eta_1(a)]$, and hence, it is possible to use the
transformation \eqref{eq:chv} to rewrite equation \eqref{eq:mainq}
in the form \eqref{eq:transeq}. Note that, in view of
\eqref{eq:re1} and \eqref{eq:re2}, the initial conditions
\eqref{eq:initu} are preserved after this transformation.
Therefore, the solution $\tilde u\equiv u(x(t))$ of
\eqref{eq:transeq} satisfies these initial conditions. According
to Lemma~\ref{le:LN}, we obtain
$$ \tilde u>0,\quad \dot{\tilde u}>0,\quad {(\tilde r\tilde h^{3}\ddot {\tilde u})}^{.}>0,
\quad\mbox{on}\quad(a, \hat\mu_1].$$ Again from
\eqref{eq:re1}-\eqref{eq:re2}, \eqref{eq:oscil1} follows. As shown
in the proof of Theorem \ref{le:oscill}, if \eqref{eq:quadform}
holds on $(a, \hat\mu_1]$, then $h'(x)>0$ on $(a, \hat\mu_1]$.
Therefore, from the second relation in \eqref{eq:re1} we get
$u''>0$ on $(a, \hat\mu_1]$. \\
By similar arguments we prove the same results for $v$.
\end{Proof}

\begin{Proof}{\bf of Theorem \ref{le:oscill}}\\
We introduce the ratios
$$
\delta_0=\frac{u}{v},\quad \delta_1=\frac{u'}{v'},\quad
\delta_2=\frac{Tu}{Tv},
$$
together with their derivatives
\begin{equation}\label{eq:ratios}
\delta'_0=-\frac{\hat\sigma}{v^2},\quad
\delta'_1=-\frac{\hat\tau}{r(v')^2},\quad
\delta'_2=\frac{p\hat\tau}{(Tv)^2}.
\end{equation}
Let
\begin{equation}\label{eq:eigen}
y_{\hat\mu_1}= u - \delta_1({\hat\mu_1})v.
\end{equation}
By Lemma \ref{le:fund}, $\delta_1({\hat\mu_1})$ is well defined.
In this case, we have $y'_{\hat\mu_1}(\hat\mu_1)=0$ and
$Ty_{\hat\mu_1}(\hat\mu_1)=\hat\tau'(\hat\mu_1)=0$. Therefore,
$y_{\hat\mu_1}$ is an eigenfunction of the boundary problem
\eqref{eq:mainq}-\eqref{eq:bc2} defined on the interval $[a,
\hat\mu_1]$. From the definition of $ \hat\mu_1$ and the initial
conditions \eqref{eq:initc1}, it follows that $\hat\tau>0$ on $(a,
\hat\mu_1]$. Thus, $\delta'_1<0$ on this interval, and hence,
$y'_{\hat\mu_1}(x)\ne0$ on $(a, \hat\mu_1(a)$. From the initial
condition $y'_{\hat\mu_1}(a)=u'(a)=1$, it follows that
$y'_{\hat\mu_1}(x)>0$ and $y_{\hat\mu_1}>0$. On the other hand,
since $Ty_{\hat\mu_1}(a)=-1$, $T'y_{\hat\mu_1}(x)>0$ on $(a,
\hat\mu_1(a)]$ and  $Ty_{\hat\mu_1}(\hat\mu_1(a))=0$, then
$Ty_{\hat\mu_1}(x)<0$ on $[a, \hat\mu_1(a))$.\\
The relations \eqref{eq:re1} and \eqref{eq:re2} yield

$$
\tilde y_{\hat\mu_1}>0,\quad \dot{\tilde y}_{\hat\mu_1}>0,\quad
 {(\tilde r\tilde
h^{3}\ddot {\tilde y}_{\hat\mu_1})}^{.}<0
$$
on $(0, \tilde\mu_1(0))$, where
$\tilde\mu_1(0)=\int_a^{\hat\mu_1}h$. From this and $\ddot {\tilde
y}_{\hat\mu_1}(0)=0$, it follows that $\ddot {\tilde
y}_{\hat\mu_1}(t)<0$ on $(0, \tilde\mu_1(0))$. It is easily seen
that, if $q\leq0$ on $[a, \hat\mu_1(a))$, then $h'<0$ on $(a,
\hat\mu_1(a))$. Therefore, from the second relation of
\eqref{eq:re1}, we obtain $y''_{\hat\mu_1}(x)<0$ on $(a,
\hat\mu_1(a))$. The theorem is proved.

\end{Proof}\\

\section {Sufficient conditions for the existence of $\hat\eta_1$}

We say equation \eqref{eq:mainq} is systems-conjugate in
$(a,\infty)$ if $\hat\eta_1$ exists; otherwise \eqref{eq:mainq} is
said to be systems-disconjugate. In this section, a number of
conjugacy and disconjugacy criteria for \eqref{eq:mainq} will be
established.

\begin{Theo}\label{le:conj1}
If $\int^\infty q(t)=-\infty$ and $\int^\infty p(t)=+\infty$ then
 equation \eqref{eq:mainq} is systems-conjugate.
\end{Theo}
\begin{Proof}
If the subwronskian $\hat\sigma$ has a zero in $(a, \infty)$, then
by Rolle's theorem, $\hat\eta_1$ exists. Assume that $\hat\sigma>0$
on $(a, \infty)$ and let $k(x)=\frac{\hat\tau'}{\hat\sigma}$. By
using the identity \eqref{eq:ident1}, we obtain
$$
k'(x) = - P -\frac{{k}^2}{r} <0 \quad\mbox{on}\qquad(a, \infty).
$$
Integrating this expression, and taking into account the
assumption $\int^{\infty}p=+\infty$, it follows that
$k(x)\rightarrow-\infty$ as $x\rightarrow+\infty$, and hence,
$\hat\mu_1(a)$ exists. Therefore, in view of Theorem \ref{le:main}
and the assumption $\int^{\infty}q=-\infty$, $\hat\eta_1(a)$
exists, which implies that \eqref{eq:mainq} is systems-conjugate.
\end{Proof}

\begin{Theo}\label{le:conj2}
If $\int^\infty \frac{1}{r}(t)=+\infty$ and $\int^\infty
q(t)=-\infty$ then equation \eqref{eq:mainq} is systems-conjugate.
\end{Theo}

For the proof of this theorem we need the following result.\\

\begin{Theo}\label{le:lei-win} {\rm (\cite{Leigh1,Win})}
If the conditions
$$\int_{a}^{\infty}r^{-1}(x)dx=\infty,\qquad\int_a^{\infty}{q(s)}ds=-\infty$$
hold, then the second-order equation \eqref{eq:disf} is oscillatory
on $(a,\infty)$; i.e., each of its solution has infinitely many
zeros in this interval.
\end{Theo}

\begin{Proof}

It is easy to see that the zeros of the subwronskian $\hat\tau'$
related to \eqref{eq:mainq} for $p\equiv0$ coincide with those of
the solution $h$ of the second-order initial value problem
\eqref{eq:disf}-\eqref{eq:initial}. In view of
Theorem~\ref{le:lei-win}, $h$ has infinitely many zeros in
$(a,\infty)$. Therefore, the first-systems focal point
$\hat\mu_1(a)$ exists for \eqref{eq:mainq} with $p\equiv0$. By
 Theorem~\ref{le:comp}, $\hat\mu_1(a)$ exists for $p>0$, and hence the assumption
$\int^{\infty}q=-\infty$ and Theorem~\ref{le:main} yield the
existence of the first-systems conjugate point $\hat\eta_1(a)$.
The theorem is proved.
\end{Proof}

By combining Theorem \ref{le:quadrat} with the second statement of
Theorem~\ref{le:main} we obtain the following criterion giving
the relation between the systems-disconjugacy of \eqref{eq:mainq}
and the sign of the associated quadratic functional.

\begin{Theo}\label{le:quadrat1}
If $\int^{\infty}q=-\infty$, then equation \eqref{eq:mainq} is
systems-disconjugate if, and only if,
$$
I[w, a, b]={\int_a^b r(w'')^2 + q(w')^2 - pw^2 dx}>0
$$
for each $b\in\(a, +\infty\)$ and each nontrivial admissible
function $w(x)$ on $[a, b]$ (i.e., $w(x)\in C^{1}[a, b]$, $w'$ is
absolutely continuous and $w''\in L_{2}[a, b]$) for which
$w(a)=w'(b)=0$.

\end{Theo}


\begin{thebibliography}{99}
\bibitem{Atk} F. V. Atkinson, {\it Discrete and Continuous Boundary Problems},
Academic Press, Newyork London 1964.


\bibitem{BK} D.Banks and G.Kurowski, A Pr\"ufer transformation for the equation
of a vibrating beam, \emph{Trans. Amer. Math. Soc.} {\bf 199}
(1974), 203--222.

\bibitem{barr} J. H. Barrett, Systems-disconjugacy of a fourth order
differential equation, \emph{Proc. Amer. Math. Soc.} {\bf 12}
(1961), 205--213.

\bibitem{barrAd} J. H. Barrett, Oscillation theory of ordinary linear differential equations
, \emph{Advan. Math.} {\bf 3} (1969), 415--509.


\bibitem{barrd} J. H. Barrett, Two point boundary Problems for self-adjoint linear differential equatios of the fourth
order with middle term, \emph{Duke Math. J. } {\bf 29} (1962),
543--554.

\bibitem{cheng} Sui-Sun Cheng, Systems-conjugate and focal points of fourth order
non-selfadjoint differential equations, \emph{Trans. Amer. Math.
Soc.} {\bf 223} (1976), 155--165.

\bibitem{cole} W. J. Coles, A general Wirtinger-type inequality, \emph{Duke Math. J. }
{\bf 27} (1960), 133--138.

\bibitem{green} L. Greenberg, An oscillation method for
fourth-order self-adjoint two point boundary value Problems with
non linear eigenvalues, \emph{Siams J. Math. Anal. } {\bf 22}
(1991), 1021--1042.

\bibitem{Leigh1} W. Leighton, On self-adjoint differential equations
of second-order, \emph{ J. London Math. Society}, {\bf 35} (1952),
37--47.

\bibitem{LN} W. Leighton, Z. Nehari, On the oscillation of solutions
of self-adjoint linear differential equations of fourth-order,
 \emph{Trans. Amer. Math. Soc.} {\bf 98} (1958), 325--377.

\bibitem{Morse} M. Morse, A generalisation of the Sturm separation and comparison theorems in n-space,
 \emph{Math. Annal.} {\bf 108} (1930), 53--69.

\bibitem{muller2} M. Pfeiffer, Oscillation criteria for self-adjoint fourth-order differntial equation,
 \emph{J. Differential Equations}, 46, 1982, p. 194--215.

\bibitem{courant} H. F. Weinberger, {\it Variational Methods for Eigenvalue Approximation}, SIAM Philadelphia, 1974.

\bibitem{Win} A. Winter, A criterion of oscillatory stability,
\emph{Quart. J. Math.},  V 7, (1949), 115-117.
\end{thebibliography}
\end{document}